\documentclass[11pt,a4paper]{article}
\usepackage[utf8]{inputenc}
\usepackage{amsmath, amssymb, amsthm}
\usepackage{graphicx,hyperref,geometry}
\usepackage{a4wide}
\usepackage{xcolor,verbatim,wasysym}
\usepackage{makecell}

\newtheorem{theorem}[subsection]{Theorem}

\newtheorem{lemma}[subsection]{Lemma}

\newtheorem{fact}[subsection]{Fact}
\newtheorem{remark}[subsection]{Remark}
\newtheorem{conjecture}[subsection]{Conjecture}

\newtheorem{question}[subsection]{Question}

\newcommand{\Z}{\mathbb{Z}_2\times\mathbb{Z}_2\times\mathbb{Z}_2}
\newcommand{\cO}{\mathcal{O}}
\newcommand{\cut}[1]{}

\title{Improved upper bound on the Frank number of $3$-edge-connected graphs}

\author{János Barát\thanks{Research supported by ERC Advanced Grant ”GeoScape” and the National Research, Development and Innovation Office, grant K-131529.} \\
\small Alfr\'ed R\'enyi Institute of Mathematics\\
\small University of Pannonia, Department of Mathematics\\
\small 8200 Veszprém, Egyetem utca 10., Hungary\\
\small \url{barat@renyi.hu} \\
and\\
Zoltán L. Blázsik\thanks{\protect\includegraphics[height=1cm]{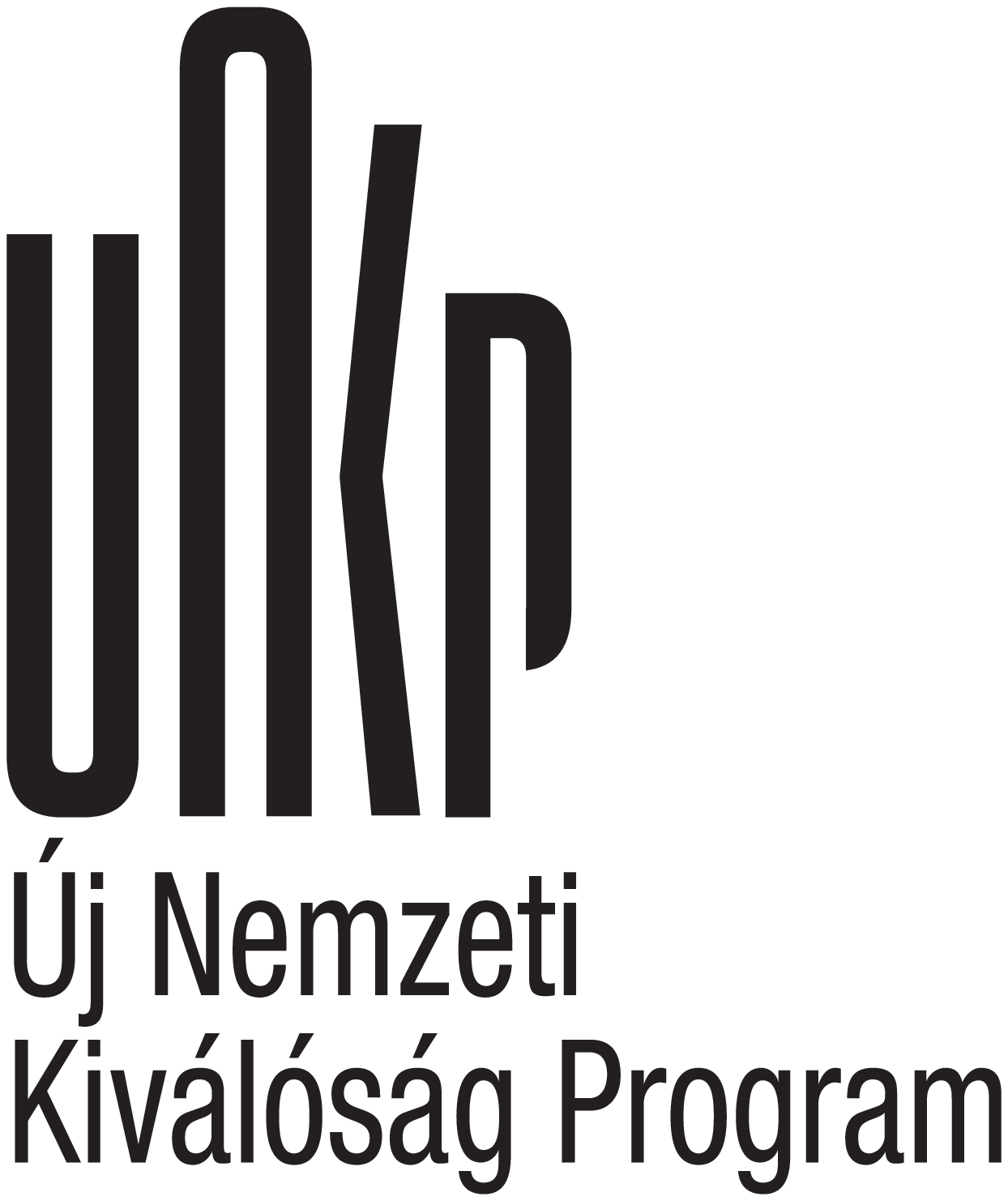}\includegraphics[height=0.8cm]{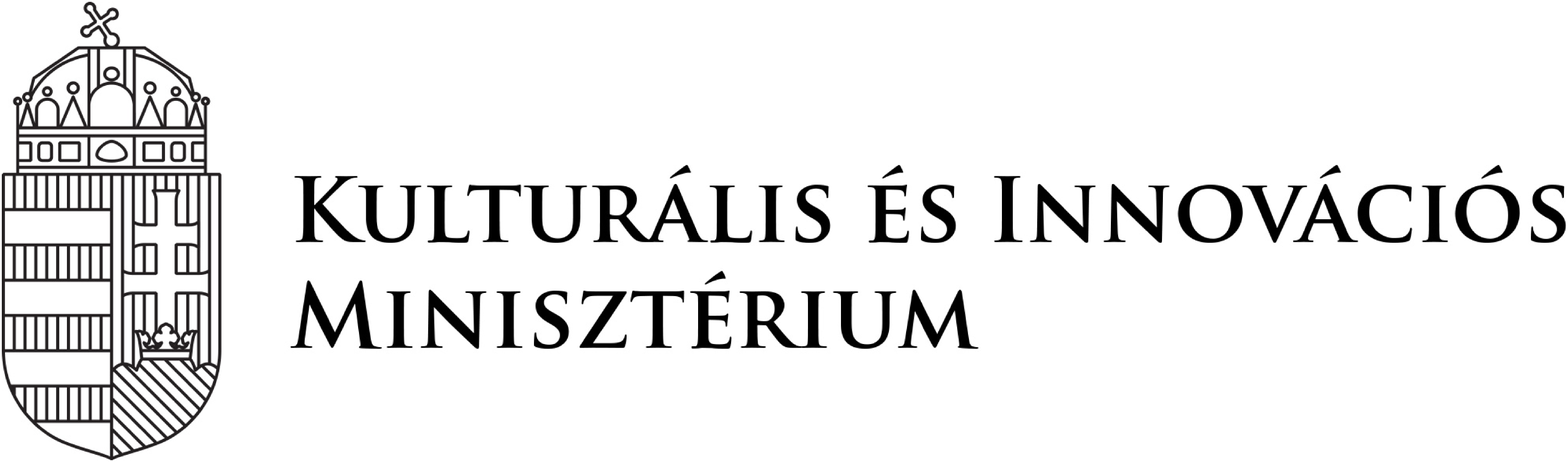} The author was supported by the \'UNKP-22-4-SZTE-480 New National Excellence Program of the Ministry for Culture and Innovation from the source of the National Research, Development and Innovation Fund. The research was supported by the Hungarian National Research, Development and Innovation Office, OTKA grant no. SNN 132625.} \\ 
\small Alfr\'ed R\'enyi Institute of Mathematics\\
\small MTA--ELTE Geometric and Algebraic Combinatorics Research Group \\
\small SZTE Bolyai Institute \\
\small \url{blazsik@renyi.hu} \\
}
\date{\today}

\begin{document}

\maketitle

\begin{abstract}
In an orientation $O$ of the graph $G$, an arc $e$ is deletable if and only if $O-e$ is strongly connected.
For a $3$-edge-connected graph $G$, the Frank number is the minimum $k$ for which $G$ admits $k$ strongly connected orientations such that for every edge $e$ of $G$ the corresponding arc is deletable in at least one of the $k$ orientations. 
Hörsch and Szigeti conjectured the Frank number is at most $3$ for every $3$-edge-connected graph $G$.
We prove an upper bound of $5$, which improves the previous bound of $7$. 
\end{abstract}

\section{Introduction}
The graphs in this paper are finite and without loops or multiple edges. We recommend the excellent book by Bondy and Murty \cite{b&m} for the concepts and notations used here.

A graph $G$ is defined by its vertex set $V$ and edge set $E$. 
An {\it orientation} of $G$ is a directed graph $D= (V,A)$ such that each edge $uv\in E$ is replaced by exactly one of the arcs $(u,v)$ or $(v,u)$. A {\it circuit} is a directed cycle. A graph is {\it cubic} if every vertex has degree $3$.
A graph is $k${\it -edge-connected} if and only if the removal of any $k-1$ edges leaves a connected graph. A $2$-edge-connected graph is often called {\it bridgeless}.

A directed graph is {\it strongly connected} if and only if selecting two arbitrary vertices $x$ and $y$, there is a directed $(x,y)$-path. 
An orientation of $G$ is {\it k-arc-connected} if and only if the removal of any $k-1$ arcs leaves a strongly connected directed graph.

\begin{theorem}[Robbins]
A graph has a strongly connected orientation if and only if it is $2$-edge-connected.
\end{theorem}

The following theorem is a fundamental result in the theory of directed graphs \cite{n-w}. 

\begin{theorem}[Nash-Williams] \label{t:NW}
A graph has a $k$-arc-connected orientation if and only if it is $2k$-edge-connected.
\end{theorem}

These results served as a motivation for H\"orsch and Szigeti \cite{szi} to investigate the following concepts. In an orientation $O$ of $G$, the arc $e$ is {\it deletable} if and only if $O{-}e$ is strongly connected. For a $3$-edge-connected graph $G$, the {\it Frank number}, denoted by $F(G)$, is the minimum $k$ for which $G$ admits $k$ strongly connected orientations such that for every edge $e$ of $G$ the corresponding arc is deletable in at least one of the $k$ orientations. Why $3$-edge-connected graphs? 
Suppose $G$ has a cut of size at most $2$. 
In any orientation of $G$, the removal of any of these edges results in either a directed cut or a directed graph that is not even connected. Hence no set of orientations can satisfy the conditions. 
On the other hand, if $G$ is $4$-edge-connected, then $G$ admits a $2$-arc-connected orientation by Theorem \ref{t:NW}. This orientation yields that $F(G)=1$ by definition. 
Consequently, $F(G)=1$ if and only if $G$ is $4$-edge-connected. Thus the problem is interesting only if the edge-connectivity of $G$ is 3, i.e. $\kappa'(G)=3$. In the sequel, we consider only graphs with edge-connectivity $3$.

H\"orsch and Szigeti \cite{szi} showed that any $3$-edge-connected graph $G$ satisfies $F(G)\le 7$.
Prior to that, DeVos et al. proved a more general result with a weaker bound \cite{unpub}:
For every 3-edge-connected graph $G$, there exists a partition of $E(G)$
into at most nine sets $\{X_1, X_2,\dots, X_m\}$ so that $G \setminus X_i$ is 2-edge-connected for every $1 \le i \le m$. Our main result improves the best known upper bound on the Frank number.

\begin{theorem} \label{fotetel}
For every $3$-edge-connected graph the Frank number is at most $5$.
\end{theorem}

The paper is organised as follows. In the second section, we introduce the main tools and results, which we use in our proof. The third 
section is dedicated to the proof
of our main result. We conclude by discussing the limits of our proof technique.

\section{Preliminaries} \label{sec:pre}

For some integer $k$, a $k${\it-flow} $(o, v)$ on a graph $G$ consists of an orientation $o$ of the edges of $G$ and a valuation $v : E(G) \mapsto \{0,\pm 1,\pm 2,\dots,\pm (k-1)\}$ such that
at every vertex the sum of the values on incoming edges equals the sum on the
outgoing edges. A $k$-flow $(o, v)$ is {\it nowhere-zero} if the value of $v$ is not $0$ for any edge of $G$. A nowhere-zero $k$-flow on $G$ is {\em all-positive} if the value $v(e)$ is positive for every edge $e$ of $G$. 
Every nowhere-zero $k$-flow can be
transformed to an all-positive nowhere-zero $k$-flow by changing the orientation of the edges with negative $v(e)$ and changing negative values of $v(e)$ to $-v(e)$. Inspired by their ideas and approach, we use the following result by Goedgebeur et al. \cite{gmr}.

\begin{lemma}\label{l:1del}
Let $G$ be a $3$-edge-connected graph, and let $(o,v)$ be an all-positive nowhere-zero $k$-flow on $G$. Any edge of $G$, which receive value $1$ in $(o, v)$ is deletable in $o$.
\end{lemma}

Indeed, an orientation arising from a flow is always strongly connected, and the removal of any arc of value 1 cannot create a directed cut since the flow is nowhere-zero. Using a slightly stronger Lemma, Goedgebeur et al. \cite{gmr} proved the following result, which we state without proof.

\begin{theorem}\label{4ff2}
If a graph $G$ admits a nowhere-zero $4$-flow, then $F(G)=2$.
\end{theorem}

It is well-known that every $3$-edge-colorable cubic graph admits a nowhere-zero $4$-flow. Tutte posed the following stronger claim, which inspired a vast amount of research.

\begin{conjecture}[Tutte's $4$-flow conjecture]\label{Tutte}
Every bridgeless graph without a Petersen-minor has a nowhere-zero $4$-flow.
\end{conjecture}

Everyone believes the validity of this conjecture.
This partly explains why the only known examples of graphs with Frank number 3 are created from the Petersen graph using certain operations \cite{frank1}. Each of the constructed graphs contains the Petersen graph as a minor. 

Since we consider 3-edge-connected graphs only, we can use the following result by Jaeger \cite{8flow}, which was later improved by Seymour \cite{6flow}. 

\begin{theorem}[Jaeger]\label{t:8flow}
Every bridgeless graph has a nowhere-zero $8$-flow.
\end{theorem}

\begin{theorem}[Seymour]\label{t:6flow}
Every bridgeless graph has a nowhere-zero $6$-flow.
\end{theorem}

Let $H$ denote an Abelian group.
An $H${\it -flow} on an oriented graph $D$ is an assignment of values of $H$ to the arcs of $D$
such that for each vertex $v$, the sum of the values on the incoming arcs is the same as the sum of the values on the outgoing arcs.
For a graph $G$, an $H$-flow is defined accordingly.
A {\it nowhere-zero} $H$-flow on $G$ is an $H$-flow, where $0 \in H$ is not
assigned to any edge.
The following is a useful corollary of a theorem by Tutte:

\begin{fact}\label{f:Abel}
If $H$ and $H'$ are two finite Abelian groups of the same order, then the graph $G$ has an $H$-flow if and only if $G$ has an $H'$-flow.
\end{fact}

In what follows, we particularly use a $\Z$-flow on a $3$-edge-connected graph $G$. 
We use $0/1$ vectors with three coordinates to denote elements of $\Z$.
The flow condition at a vertex $v$ implies the following property: in each coordinate, the sum of values on edges incident to $v$ is $0$. Hence there are an even number of $1$'s in each coordinate, on the edges incident to a fixed vertex.

\section{Improvement of the upper bound}

Surprisingly, Theorem~\ref{t:8flow}, the weaker flow result of Jaeger is the useful for our purposes.
The main idea of the proof is the following. 
We fix a nowhere-zero 8-flow of a $3$-edge-connected graph $G$ existing by Theorem \ref{t:8flow}.
We create 5 other nowhere-zero $k$-flows of $G$ in such a way that we control the set of edges with value 1, and we apply Lemma~\ref{l:1del}. Since every edge of $G$ receives value 1 in at least one of the flows, we are done. Let us recall our main theorem.

\medskip

\noindent{\bf Theorem~\ref{fotetel}}
For every $3$-edge-connected graph the Frank number is at most $5$.

\begin{proof}
Combining Theorem \ref{t:8flow} and Fact \ref{f:Abel}, we consider a nowhere-zero $\Z$-flow on $G$. For $i\in\{1,2,3\}$, let $G_i$ denote the subgraph of $G$ induced by those edges of $E$, which have value 1 at the $i$th coordinate. Note that these subgraphs might have some edges in common, more precisely the number of nonzero coordinates of each edge is the same as the number of subgraphs to which it belongs. By the nowhere-zero property, every edge is contained in at least one of these subgraphs.

By the flow condition, the subgraphs $G_1,G_2,G_3$ are Eulerian. We can think of an Eulerian trail as a directed graph. 
Thus we can partition each edge set $E(G_i)$ into edge-disjoint circuits.
We may use some vertices in different Eulerian trails. 
For $i=1,2,3$, let us fix an orientation $o_i$ of $E(G_i)$. 
At this point, it might occur that an edge of $G$ has different orientations in different subgraphs $G_i$. 
The next aim is to select, for each $i\in\{1,2,3\}$, an appropriate positive value $v_i$ to send along $o_i$. 

Let us emphasize that after fixing the pair $(o_i,v_i)$ for each $G_i$, we define another flow $(O_1,f_1)$ in the next step.
The orientation and value of an arc in $(O_1,f_1)$ is determined by the 
superposition of the chosen $(o_i,v_i)$ pairs or triples.
We always assign a positive value. 
If the orientations go opposite, then we let the largest value determine the direction and subtract the values going in opposite direction. 

Let us construct an all-positive nowhere-zero 8-flow $(O_1,f_1)$ by sending values $v_1=1$, \mbox{$v_2=2$} and $v_3=4$ along the fixed orientations $o_i$ of $G_i$, respectively. Indeed, this is a flow, and there can be no arcs of value 0. 
The maximum value of an arc is 7, if this edge was of type $(1,1,1)$ in the original $\Z$-flow, and the edge received the same orientation in all three orientations $o_i$. 
Since in every edge-cut, the sum of the flow values in the two directions is the same, there are no directed cuts if we take the superposition of the flows $o_1,o_2,o_3$. Thus $O_1$ is strongly connected.

In order to prove $F(G)\le 5$, we create a set of strongly-connected orientations $\cO$ of $G$ such that $|\cO|=5$, and any edge of $G$ is deletable in at least one of the orientations of $\cO$. Let $O_1$ be the first member of $\cO$. By Lemma \ref{l:1del}, the arcs of value 1 in $(O_1,f_1)$ are deletable with respect to $O_1$. What are these arcs? We summarize the possible final orientations and values of the arcs in Table \ref{tab:type}.

\begin{table}[!h]
    \centering
    \begin{tabular}{|c|c|}
        \hline
        \makecell{Value of $e$ \\ in the original \\ $\Z$-flow} & Final orientation and value of $e$ in $(O_1,f_1)$ \\ \hline
        $(1,0,0)$   & same as in $o_1$, 1 \\ \hline
        $(0,1,0)$   & same as in $o_2$, 2 \\ \hline
        $(0,0,1)$   & same as in $o_3$, 4 \\ \hline
        $(1,1,0)$   & $\left \{ \begin{tabular}{cc}
            same as in $o_2$, 3 & \qquad if $e$ has the same orientation in $o_1$ and $o_2$, \\
            same as in $o_2$, 1 & \qquad if $e$ has different orientations in $o_1$ and $o_2$.             
        \end{tabular} \right .$ \\ \hline
        $(1,0,1)$   & $\left \{ \begin{tabular}{cc}
            same as in $o_3$, 5 & \qquad if $e$ has the same orientation in $o_1$ and $o_3$, \\
            same as in $o_3$, 3 & \qquad if $e$ has different orientations in $o_1$ and $o_3$.             
        \end{tabular} \right .$ \\ \hline
        $(0,1,1)$   & $\left \{ \begin{tabular}{cc}
            same as in $o_3$, 6 & \qquad if $e$ has the same orientation in $o_2$ and $o_3$, \\
            same as in $o_3$, 2 & \qquad if $e$ has different orientations in $o_2$ and $o_3$.
        \end{tabular} \right .$ \\ \hline
        $(1,1,1)$   & $\left \{ \begin{tabular}{cc}
        same as in $o_3$, 7 & \qquad \makecell{if $e$ has the same orientation \\ in all three of $o_1,o_2,o_3$,} \\
            same as in $o_3$, 5 & \qquad \makecell{if the orientation of $e$ in $o_1$ is different\\ from that in $o_2$ and $o_3$,} \\
            same as in $o_3$, 3 & \qquad \makecell{if the orientation of $e$ in $o_2$ is different\\ from that in $o_1$ and $o_3$,} \\
            same as in $o_3$, 1 & \qquad \makecell{if the orientation of $e$ in $o_3$ is different\\ from that in $o_1$ and $o_2$.}
        \end{tabular} \right .$ \\ \hline
    \end{tabular}
    \caption{The possible orientations and values of $e$ in $(O_1,f_1)$.}
    \label{tab:type}
\end{table}

There are three types of arcs of value 1 in $(O_1,f_1)$. 
We create 4 other flows such that the larger valued arcs of $O_1$ receive value 1 in at least one orientation. We refer to any edge of $G$ as one which belongs to the corresponding row of Table \ref{tab:type}.
The next goal is to create 4 other strongly connected orientations such that every type of edge defined by Table \ref{tab:type} gets flow value 1 in at least one of those orientations, and hence is deletable in those orientations by Lemma~\ref{l:1del}. 
We achieve this by combining the following two things: we reverse the orientations $o_i$ in $G_i$ for each $i\in\{1,2,3\}$ if necessary, and we change the values $v_i$ appropriately. 
By reversing the orientations, we can switch the role of a fixed edge with respect to its role described in Table \ref{tab:type}, and by changing the values we can change the final values of the edges.

We define the following two $10$-flows ($O_2$,$O_3$) and two $8$-flows ($O_4$,$O_5$) based on the orientations $o_1,o_2,o_3$. Again, if an edge is present in several circuits, we superimpose the values as before\footnote{the larger value deciding the direction except if 2+3 goes in one direction and 4 in the opposite} to make the flow $(O_i,f_i)$ all-positive for all $i\in\{1,2,3,4,5\}$. We define $O_2$ by keeping the same orientations as in $O_1$, but sending value 4 along $o_1$, value 2 along $o_2$, and value 3 along $o_3$. The only difference between $O_3$ and $O_2$ is that the orientations of the circuits in $o_3$ are reversed. We get $O_4$ from $O_1$ by reversing the orientation of the circuits in $o_2$ and sending value 2 along $o_1$, value 1 along $o_2$ and value 4 along $o_3$. The only difference between $O_5$ and $O_4$ is that the values along $o_2$ and $o_3$ are swapped. We indicate in Table~\ref{tab:remain} which edges receive value 1 in the corresponding orientations of those edges, which did not receive value 1 with respect to $O_1$.

\begin{table}[!h]
    \centering
    \begin{tabular}{|c|c|c|c|c|}
    \hline
    original value & $(O_2,f_2)$ & $(O_3,f_3)$ & $(O_4,f_4)$ & $(O_5,f_5)$ \\
    \hline
    $(0,1,0)$ & & & value 1 & \\
    \hline
    $(0,0,1)$ & & & & value 1 \\
     \hline
    $(1,1,0)$ & & & value 1& \\
    \hline
    $(1,0,1)$ row 6 & & value 1& & \\
     \hline
    $(1,0,1)$ row 7 & value 1 & & & \\
    \hline
    $(0,1,1)$ row 8 & & value 1 & & \\
       \hline
    $(0,1,1)$ row 9 & value 1 & & &\\
     \hline
    $(1,1,1)$ row 10 & & & & value 1\\
     \hline
     $(1,1,1)$ row 11 & value 1 & & & \\
       \hline
    $(1,1,1)$ row 12 & & value 1 & & \\
     \hline
   \end{tabular}
    \caption{Cross-table for clearing the remaining cases referring to Table~\ref{tab:type}}
    \label{tab:remain}
\end{table}

By Lemma \ref{l:1del}, it is clear that $\cO=\{O_1,O_2,O_3,O_4,O_5\}$ yields $F(G)\le 5$. 
\end{proof}

Let us continue with some comments and observations on the proof.  

\begin{remark}\label{r:flow}
It does not really matter what kind of flow we use to construct an orientation of $\cO$. However, there are two key observations: since the orientation arises from flows, it cannot have a directed cut. Therefore, it is a strongly connected orientation, and by the nowhere-zero property and the $3$-edge-connectivity of $G$ an arc of value 1 cannot be the only arc going in the opposite direction in any cut. 
\end{remark} 

\begin{remark}
Note that we do not claim that these are the only deletable edges of $O_i\in\cO$. It may happen that some arcs with higher values are also deletable. Hence our result is slightly stronger than just proving $F(G)\le 5$ since we only relied on the arcs of value 1. 
\end{remark}

\section*{Discussion}

We do not see how to use directly Theorem~\ref{t:6flow} to improve the Frank number of $3$-edge-connected graphs. A natural question arises, can we use a similar proof technique to achieve an even better upper bound? Goedgebeur et al. \cite{gmr} proved that if $G$ admits a nowhere-zero all-positive 4-flow, then $F(G)\le 2$. 
On the other hand, if we start with an arbitrary all-positive $k$-flow for some $k\ge5$, then we either need to have some control over the arcs with value 1 or need to know some underlying structure, which we can exploit. 
In this proof technique, it seems vital to have those Eulerian subgraphs. To find them, we had to use an equivalent $H$-flow for an appropriate group $H$. However, the even degrees within those subgraphs was guaranteed by the fact that the corresponding coordinate was 0/1 valued. Hence it is unclear to us, what kind of $H$ can be used here other than $\mathbb{Z}^m_2$.

Suppose that $H=\mathbb{Z}^m_2$ for some $m\ge2$. Let us introduce the notation $A_k$ for those arcs of the original $H$-flow, which have exactly $k$ non-zero coordinates ($k\in\{1,2,\dots,m\}$). Similarly to our proof, an all-positive $H$-flow gives rise to $m$ Eulerian subgraphs $G_i$ (for $i\in\{1,2,\dots,m\}$) and we can fix an orientation $o_i$ on each of them. How many different types of arcs can get value 1? 
We carefully choose values $v_i$ and consider orientations $(O_j,f_j)$ as the superposition, while no arc gets value 0 in the process.

The existence of the arcs of $A_2$ ensures that the values $v_i$ must be pairwise different, otherwise there is a possibility that an arc gets value 0 in the superposition. 
Thus at most one of the families in $A_1$ can get value 1 (the ones within $G_i$ if and only if $v_i=1$) in each orientation. Hence with this technique, we definitely need at least $m$ orientations, consequently $m\le F(G)$.

There are $2\binom{m}{2}$ possible types of arcs corresponding to the arcs of $A_2$.
At most two of them can have value 1 in a fixed orientation $(O,f)$. 
Indeed, if there were at least three types of such arcs, then we would choose two such that the first one is non-zero at coordinates $i_1,i_2$, the second one at coordinates $j_1,j_2$ such that $v_{i_1}-v_{i_2}=1=v_{j_1}-v_{j_2}$, and $i_1,i_2,j_1,j_2$ are pairwise different. This could possibly lead to some arcs with final value 0, whenever an arc $xy$ of $A_4$ has $1$'s at coordinates $i_1,i_2,j_1,j_2$. In addition, the corresponding four orientations are chosen in such a way that $xy$ is oriented in the same direction in $G_{i_1}$ and $G_{j_2}$, but in the opposite direction in the other two Eulerian subgraphs. 
That yields $\binom{m}{2}\le F(G)$. 
Thus if one would like to improve our result, then $m$ should be at most 3. Since $m=2$ corresponds to the $4$-flows, that circles back to Conjecture \ref{Tutte}, the only possibility is $m=3$.

From the previous observations, one can deduce that we need at least three orientations to deal with the arcs of $A_1$. 
In each of these orientations, there is an Eulerian subgraph $G_i$ for $i\in\{1,2,3\}$ such that $v_i=1$. 
These orientations cannot have two types of arcs with value 1 from $A_2$, since that would be possible only if the three values were $1$, $2$ and $3$. Again this contradicts the non-zero property of our flow ($1+2-3=0$).
Consequently, at least three types of arcs are missing after these three orientations from the arcs of $A_2$. Hence at least two other orientations are needed. Therefore our upper bound of 5 is best if $H=\mathbb{Z}_2\times \mathbb{Z}_2 \times \mathbb{Z}_2$. We would like to conclude by posing a natural problem.

\begin{question}
    Is it possible to improve the upper bound on $F(G)$ using a similar flow technique but starting with a different group $H$?
\end{question}

\end{document}